\documentclass[reqno,a4paper,11pt]{amsart}
\usepackage[english]{babel}
\usepackage{graphicx,a4wide,mathrsfs,tikz,mathtools,latexsym,ifthen,amssymb,stmaryrd}
\mathtoolsset{showonlyrefs,showmanualtags}
\usepackage{enumerate}
\usepackage{times}

\usetikzlibrary{positioning,shapes,shadows,arrows}
\usepackage[utf8x]{inputenc}
\usepackage{lmodern} 
\usepackage[toc,page]{appendix}
\usepackage{todonotes}

\usepackage{color}
\definecolor{MyDarkBlue}{rgb}{0.15,0.25,0.45}
\usepackage[pdftex]{hyperref}
\hypersetup{
colorlinks=true,
citecolor=MyDarkBlue,
linkcolor=black,
urlcolor=black,
pdfauthor={},
pdftitle={},
pdfsubject={}
breaklinks=true
}

\newcommand{\triend}{\parbox{2mm}{\hfill} \hfill\mbox{\hspace{0.2mm}}\hfill$\triangle$}
\newcommand{\ocend}{\parbox{2mm}{\hfill} \hfill\mbox{\hspace{0.2mm}}\hfill$\oslash$}

\newtheorem{theorem}{Theorem}[section]
\newtheorem{proposition}[theorem]{Proposition}
\newtheorem{lemma}[theorem]{Lemma}
\newtheorem{corollary}[theorem]{Corollary}
\newtheorem*{theorem*}{Theorem}

\theoremstyle{remark}
\newtheorem{example}[theorem]{Example}

\theoremstyle{remark}
\newtheorem{rem}[theorem]{Remark}
\newenvironment{remark}{\begin{rem}}{\triend\end{rem}}

\theoremstyle{definition}
\newtheorem{defin}[theorem]{Definition}
\newenvironment{definition}{\begin{defin}}{\ocend\end{defin}}

\title{Instantons and framed sheaves on K\"ahler Deligne-Mumford stacks}

\author{Philippe Eyssidieux, Francesco Sala} 

\address{Universit\'e Joseph Fourier et Institut Universitaire de France, France}

\email{Philippe.Eyssidieux@ujf-grenoble.fr}
  
\address{Institut des Hautes \'Etudes Scientifiques, France}

\email{salafra83@gmail.com}

\date{\today}
\begin{document}

\maketitle

\begin{abstract}
We provide stacky generalizations of classical gauge-theoretic results inspired by Donaldson, the Uhlenbeck-Yau theorem and variants due to Bando and his collaborators. Moreover, we show an application of this machinery in the study of ALE spaces.
\end{abstract}

\tableofcontents

\section{Introduction}

It is a general principle (or a {\em metatheorem}), which was explained to us by C. Simpson, that the basic results of  
Differential Geometry (Real or Complex) obtained through the analysis of Geometric Partial Differential 
Equations extend to {\em Deligne-Mumford} stacks\footnote{Complex differential geometry on {\em Artin} stacks is much more delicate 
but may have more applications - however we will not pursue this direction here.}. The main obstacle
to carry out this extension is to set up the stacky definitions correctly. The orbifold versions are (in principle) 
well-known but nontrivial generic isotropy is not a very serious technical obstacle. 
Doing the proof of such a result for stacks  amounts to checking that the extension to stacks of the constructions used in the proof is possible and sufficiently 
well documented in the literature. 
In this article, we enforce this principle and state stacky generalizations 
of classical gauge-theoretic results inspired by Donaldson, the Uhlenbeck-Yau theorem \cite{art:uhlenbeckyau1986} and variants due to Bando and his collaborators.

Our policy with respect to  writing-up these results will be to give  no details on their deep
gauge-theoretic or analytical aspects and concentrate on purely stack-theoretic issues.  Perhaps, one may regret that no 
comprehensive reference on global analytic methods on stacks in the style of \cite{book:laumonmoret-bailly2000} is available. If more applications of this method in the realm of ordinary complex/differential geometry are developed, it may become inevitable to lay down these foundations.

Our motivation was to generalize the following correspondence due to Donaldson to 
instantons on ALE spaces. In \cite{art:donaldson1984} Donaldson proved that there is a one-to-one correspondence 
between  $U(r)$-instantons of charge $n$ on $\mathbb R^4$ (modulo gauge equivalence) and {\em framed vector bundles} 
of rank $r$ and second Chern class $n$ on the complex projective plane $\mathbb{P}^2=\mathbb{C}^2\cup l_\infty$. 
A framed vector bundle is a pair $(E,\phi_E)$, where $E$ is a vector bundle on $\mathbb{P}^2$ of rank $r$
and $c_2(E)=n$ and $\phi_E\colon E_{\vert l_{\infty}}\xrightarrow{\sim}\mathcal{O}_{l_{\infty}}^{\oplus r}$ 
a trivialization along the line $l_{\infty}$. Donaldson's correspondence was generalized to  instantons on the 
blowup $\widetilde{\mathbb{C}^2}$ of $\mathbb{C}^2$ at $n$ points (for $n=1$ by King \cite{phd:king1989} and for $n\geq 2$ by Buchdahl \cite{art:buchdahl1993}). 
In all these examples, in order to state a correspondence between instantons and framed 
vector bundles one needs first to endow the noncompact manifold of a K\"ahler structure, compactify it by adding a projective line, 
and consider framed vector bundles on the corresponding smooth projective surface. 

A natural question to ask is if this algebro-geometric approach can
be applied to other noncompact 4-dimensional Riemannian manifolds. In particular, we are interested in ALE spaces of type $A_{k-1}$, where $k\geq 2$ is an integer. A first result in
the ALE case is due to Bando: in \cite{art:bando1993}, as a consequence of a generalization of Donaldson's characterization of instantons on $\mathbb R^4$ to general K\"ahler manifolds, he proved that if one can compactify an ALE space $X^{\circ}$ of type $A_{k-1}$ by adding a smooth divisor $D$ which has positive normal
line bundle so that one obtains a compact K\"ahler manifold $X:=X^{\circ}\cup D$, the following correspondence holds: 
holomorphic vector bundles on $X$ unitary flat along $D$ are in a one-to-one correspondence with holomorphic vector bundles on $X^{\circ}$ 
endowed with a Hermite-Einstein metric with trivial holonomy at infinity. If such compactification $X$ exists, the smooth curve $D$ would be connected (since $X^\circ$ has only one end),
would have genus zero (because the fundamental group of this end is finite) hence the holonomy at infinity would be the trivial representation of the 
fundamental group of the end. But this fundamental group is non trivial. The way we propose to circumvent this obstacle is to consider orbifold compactifications of $X^\circ$.

Let us fix $X^\circ=X_k$, where $X_k$ is the minimal resolution of the $A_{k-1}$ toric singularity of $\mathbb{C}^2/\mathbb{Z}_k$. In \cite{art:bruzzopedrinisalaszabo2013} a 
compactification $\mathscr{X}_k$ of $X_k$ is constructed, which turns out to be a projective toric orbifold. 
In particular $\mathscr{X}_k\setminus X_k$ is a smooth effective Cartier divisor $\mathscr{D}_\infty$ with ample and positive normal 
line bundle. 

In this paper we generalize Bando's result by proving the following:

\begin{theorem*}
There is a one-to-one correspondence
between holomorphic vector bundles on $\mathscr{X}_k$, which are isomorphic along $\mathscr{D}_\infty$ 
to a fixed vector bundle $\mathcal F_\infty$ endowed with a flat unitary connection $\nabla$,
and vector bundles on $X_k$ endowed with an  Hermite-Einstein metric such that the curvature is square integrable and the holonomy at infinity is given by the holonomy of $\nabla$.
\end{theorem*}

Let us describe the organisation of this article. In Sections \ref{sec:topstacks} and \ref{sec:diffstacks}, we set up the definitions of the basic ingredients of the differential geometry of Deligne-Mumford stacks with an emphasis on the K\"ahler case. In Sections \ref{sec:hermite-einstein} and \ref{sec:hermite-einstein2}, we extend to Deligne-Mumford stacks the Uhlenbeck-Yau and Bando theorems. In Section \ref{sec:application}, we recall the construction of the orbifold compactification of $X_{k}$ and conclude with the proof of the above theorem (Theorem \ref{thm:ALE}).

\subsection{Acknowledgements} 

We thank C.\ Simpson for useful suggestions and for his interest in our work and U.\ Bruzzo for reading and commenting on a draft of this paper. We are also grateful to M.\ Pedrini, R.\ J.\ Szabo and R.\ Thomas for helpful discussions and correspondence. The last draft of the paper was written while the second author was visiting IH\'ES and Universit\'e Joseph Fourier. He thanks both institutions for hospitality and support.

\bigskip\section{Topological, differentiable and smooth analytic stacks}\label{sec:topstacks}

In this section we briefly describe topological, differentiable and smooth analytic stacks. Our main references are 
\cite{art:behrend2004, art:noohi2005, art:behrendnoohi2006}. We assume that the reader is familiar with the notions of {\em category fibered in groupoids} and of {\em stack} (cf.\ \cite[Sect.~2 and 3]{book:laumonmoret-bailly2000}).

\subsection{Topological stacks}

Let $\mathbf{Top}$ be the category of topological spaces. We fix a final object $\ast$, \emph{the} point, in the category $\mathbf{Top}$ of topological spaces. We endow $\mathbf{Top}$ with the usual Grothendieck topology (covers are simply topological open covers); so we can talk about (the 2-category of) stacks over $\mathbf{Top}$: it is closed under fibre products and by Yoneda's lemma the category of topological spaces embeds as a full subcategory of the 2-category of stacks.

We say a morphism $f \colon \mathscr{Y} \to \mathscr{X}$ of stacks over  $\mathbf{Top}$ is {\em representable}, if for any map $U \to \mathscr{X}$ from a
topological space $U$ to $\mathscr{X}$, the fiber product $V:=U\times_{\mathscr{X}}\mathscr{Y}$ is equivalent to a topological space. 

Several properties $\mathbf{P}$ of maps of topological spaces are stable under base change (cf.\ \cite[Sect.~4.1]{art:noohi2005}). For example: \emph{to be open maps, epimorphisms, surjective maps, embeddings, closed embeddings, open embeddings, local homeomorphisms, covering maps, maps with finite fibers, maps with discrete fibers}. We say a representable morphism $f \colon \mathscr{Y} \to \mathscr{X}$ of stacks over $\mathbf{Top}$ satisfies a property $\mathbf{P}$ if for any map $U \to \mathscr{X}$ from a topological space $U$ to $\mathscr{X}$, the base extension $V \to U$ of $f$ satisfies $\mathbf{P}$. 

\begin{definition}
A {\em pre-Deligne-Mumford topological stack} is a stack $\mathscr{X}$ for which 
there exists an epimorphism $u\colon U \to \mathscr{X}$ from a topological space $U$, such that $u$ is representable by local homeomorphisms. 
We call the pair $(U,u)$ an {\em atlas} of $\mathscr{X}$. A pre-Deligne-Mumford topological 
stack $\mathscr{X}$ is called a {\em Deligne-Mumford topological stack} if $U$ is  Hausdorff  for some atlas $(U,u)$ and 
the diagonal $\mathscr{X} \to \mathscr{X}\times\mathscr{X}$ is representable by proper maps with closed discrete fibers.
\end{definition}
One can define a {\em site} $\mathcal{S}(\mathscr X)$ of a Deligne-Mumford topological stack $\mathscr X$ in the following way. The objects of the underlying category of $\mathcal{S}(\mathscr X)$ are the atlases $(U,u)$ of $\mathscr{X}$, the arrows are the morphisms $(\varphi, \alpha)\colon (U,u)\to (V,v)$ of two atlases where $\varphi\colon U\to V$ is a local homeomorphism of topological spaces and $\alpha\colon u\xrightarrow{\sim} v\circ \varphi$ is a 2-isomorphism. The topology on this category is the one induced by the pre-topology, where the {\em covering families} are of the following form: for an atlas $(U,u)$, we denote by $\mathrm{Cov}(U,u)$ the set of families of morphisms $\{(\varphi_i, \alpha_i)\colon (U_i, u_i)\to (U,u)\}_{i\in I}$ such that the map
\begin{equation}
\bigsqcup_{i\in I}\varphi_i \colon \bigsqcup_{i\in I} U_i \to U
\end{equation}
is open and surjective.
\begin{remark}
The previous definition of pre-Deligne-Mumford topological stacks comes from \cite[Sect.~3.1]{art:behrendnoohi2006}. 
In \cite[Sect.~7]{art:noohi2005}, Noohi gives a more general definition of {\em pre-topological stacks}. In \cite{art:noohi2005}, our pre-Deligne-Mumford topological 
stacks are called {\em weak Deligne-Mumford topological stacks} (cf.\ \cite[Def.~14.3]{art:noohi2005}). On the other hand our definition of Deligne-Mumford
topological stacks is more restrictive than the one given in \cite{art:behrendnoohi2006}: indeed, in \cite[Def.~3.1]{art:behrendnoohi2006} the authors do not require existence of a Hausdorff atlas 
and require only
that the diagonal is representable by closed maps with discrete finite fibers. As the authors point out in \cite[Sect.~3.1]{art:behrendnoohi2006} closedness is a property invariant under base extension only via local homeomorphism. On the other hand, we impose properness instead of closeness
because we need a property invariant under {\em any} base extension. 
\end{remark}

For a stack $\mathscr{X}$, we denote by $\pi_0\mathscr{X}$ the sheaf associated to the presheaf of sets on $\mathbf{Top}$ defined by 
$W\mapsto \{ \mbox{isomorphim classes in } \mathscr{X}(W)\}$. To any Deligne-Mumford topological stack $\mathscr{X}$ one can associate a topological space
 $X$, called the {\em coarse moduli space} of $\mathscr{X}$: as a set $X$ is equal to $\pi_0\mathscr{X}(\ast)$. For any open substack $\mathscr{X}'\subseteq \mathscr{X}$ (i.e. any representable open embedding $\mathscr{X}'\to \mathscr{X}$), we have a natural inclusion of coarse moduli spaces $X'\subseteq X$. These are defined to be the open sets of $X$.

\begin{definition}
Let $\mathscr{X}$ be a Deligne-Mumford topological stack. We say that $\mathscr{X}$ is {\em connected} if it has no proper open–closed substack. We say that $\mathscr{X}$ is \emph{compact} if its coarse moduli space $X$ is compact.
\end{definition}

The natural constructions of homotopy theory of topological spaces can be extended to Deligne-Mumford topological stacks. For instance, one can define a reasonable notion of homotopy  between maps, and this allows us to define the n-th homotopy group  of a pointed topological stack $(\mathscr{X}, x)$ as pointed homotopy  classes of maps $[(S^n, \ast), (\mathscr{X}, x)]$ (cf.\ \cite[Sect.~17]{art:noohi2005}).

As explained in \cite{art:behrend2004}, there is a well defined singular (co)homology theory with integral coefficients for Deligne-Mumford topological stacks. The (co)homology theory with integral coefficients depends on the stacky structure, unlike the (co)homology theory with rational coefficients. Indeed, one has the following result.
\begin{proposition}\cite[Prop.~36]{art:behrend2004}
Let $\mathscr X$ be a Deligne-Mumford topological stack and $X$ its coarse moduli space. Then the coarse moduli space morphism $\pi\colon \mathscr X\to X$ induces isomorphisms
\begin{equation}
H^k(\mathscr X; \mathbb{Q})\simeq H^k(X; \mathbb{Q})\ .
\end{equation}
\end{proposition}

\subsection{Analytic and differentiable stacks}\label{sec:analdiffstacks}

Let $\mathbf{Comp}$ be the category of complex manifolds\footnote{We are assuming that manifolds have a countable basis for their topology.}, endowed with the usual Grothendieck topology (where covers are simply topological open covers). 
As in the case of topological stacks, we can construct the 2-category of stacks over $\mathbf{Comp}$, and by Yoneda's lemma 
this category contains the category of complex manifolds as a full subcategory.

We say a morphism $f\colon \mathscr{Y} \to \mathscr{X}$ of stacks over $\mathbf{Comp}$ is {\em representable by local biholomorphisms}, if for any map $U
\to \mathscr{X}$ from a complex manifold $U$ to $\mathscr{X}$, the fiber product $V:=U\times_{\mathscr{X}}\mathscr{Y}$ is equivalent to a complex manifold, and the map $V \to U$ is a local biholomorphism.

\begin{definition}
A stack $\mathscr{X}$ over $\mathbf{Comp}$ is called a {\em smooth pre-Deligne-Mumford analytic stack}
if there exists an epimorphism $u \colon U \to \mathscr{X}$ from a complex manifold $U$ such that $u$ is
representable by local biholomorphisms. We call the pair $(U,u)$ an {\em atlas} of $\mathscr{X}$.
\end{definition}
\begin{rem}
Since a holomorphic map of complex manifolds is a local biholomorphism if and only if it is a local homeomorphism, 
the previous definition coincides with the one given in \cite[Sect.~3.2]{art:behrendnoohi2006}.\hfill$\triangle$
\end{rem}

We say a morphism $f \colon \mathscr{Y} \to \mathscr{X}$ of smooth pre-Deligne-Mumford analytic stacks is {\em representable}, if for any map $U \to \mathscr{X}$ from a complex manifold $U$ to $\mathscr{X}$ that is representable by local biholomorphisms, the fiber product $V:=U\times_{\mathscr{X}}\mathscr{Y}$ is equivalent to a complex manifold. Let $\mathbf{P}$ be a property of morphisms of complex manifolds that is invariant under base change with respect to local biholomorphisms. For example we can take $\mathbf{P}$=closedness, $\mathbf{P}$=to have finite fibers, $\mathbf{P}$=to have discrete fibers, $\mathbf{P}$=properness, $\mathbf{P}$=to be unramified, or $\mathbf{P}$=to be a covering space. Then, we say a representable morphism $f \colon \mathscr{Y} \to \mathscr{X}$ of smooth pre-Deligne-Mumford analytic stacks is $\mathbf{P}$, if for any map $U \to \mathscr{X}$ from a complex manifold $U$ to $\mathscr{X}$ that is representable by local biholomorphisms, the base extension $V \to U$ is $\mathbf{P}$.

\begin{definition}
A smooth pre-Deligne-Mumford analytic stack $\mathscr{X}$ is called a {\em smooth Deligne-Mumford analytic stack} 
if, for some atlas $(U,u)$ of $\mathscr{X}$, $U$ is Hausdorff and the diagonal $\mathscr{X} \to \mathscr{X}\times\mathscr{X}$ 
is representable by proper unramified\footnote{The attribute ``unramified" means ``injective differential" from the differential geometric point of view.} finite maps. An {\em orbifold} is a smooth Deligne-Mumford analytic stack with generically trivial stabilizer.
\end{definition}
Also in this case, one can define a site $\mathcal{S}(\mathscr X)$ in a way similar than before.
\begin{remark}
Our definition of smooth Deligne-Mumford analytic stacks is more restrictive that the one given in \cite{art:behrendnoohi2006}: in \cite[Def.~3.3]{art:behrendnoohi2006} the authors assume that the diagonal is representable only by closed maps with finite fibers. On the other hand, one can give a definition of Deligne-Mumford analytic stacks over the category $\mathbf{Analytic}$ 
of all analytic spaces (cf.\ \cite{art:noohi2005}). As pointed out in  \cite[Rem.~3.4]{art:behrendnoohi2006}, the 2-category of 
smooth Deligne-Mumford analytic stacks defined above is equivalent to the sub 2-category of the 2-category of Deligne-Mumford analytic stacks 
of \cite{art:noohi2005} consisting of smooth Deligne-Mumford analytic stacks. Since the Deligne-Mumford analytic stacks we are interested in are smooth, we prefered not to work over $\mathbf{Analytic}$ since this category may be less familiar to the reader.
\end{remark}

By substituting in the previous definitions the term ``local biholomorphism" with ``local diffeomorphism", a definition of smooth (pre-)Deligne-Mumford differentiable stacks over the category $\mathbf{Diff}$ of differentiable manifolds can be given and was actually introduced in \cite{art:behrend2004}. 

One can define the coarse moduli space $X$ of a smooth Deligne-Mumford analytic (resp. differentiable) stack $\mathscr{X}$. 
It is an analytic space, but it may not in general be a smooth complex manifold. 

\begin{proposition}\cite{art:behrend2004,art:noohi2005}\label{prop:localstructurestacks}
Let $\mathscr{X}$ be a Deligne-Mumford topological (resp.\ smooth differentiable, resp.\ smooth analytic) stack. 
Then there is a covering $\{\mathscr{U}_i\}$ of $\mathscr{X}$ by open substacks such that each $\mathscr{U}_i$ is a quotient stack $[Z/G]$,
where $Z$ is a Hausdorff topological space (resp. a Hausdorff differentiable manifold, resp. Hausdorff \ complex manifold), 
and $G$ a finite group acting continuously (resp. differentiably, resp.\ analytically) on $Z$.
\end{proposition}

As explained in \cite[Cor.~25]{art:behrend2004}, there is a well-defined de Rham cohomology theory for any smooth Deligne-Mumford differentiable stack $\mathscr X$. 
\begin{definition}
Let $\mathscr X$ be a smooth Deligne-Mumford differentiable stack. We say that $\mathscr X$ is {\em oriented} if for any atlas $(U,u)$ of $\mathscr X$ the associated groupoid 
\begin{equation}
  \begin{tikzpicture}[xscale=2.8,yscale=-.7, ->, bend angle=25]
\node (A0_1) at (0,1) {$U\times_{\mathscr X} U$};
\node (A1_1) at (1.15,1) {$U$}; 
\node (B1) at (0.3,0.8) {$ $};
\node (B2) at (1,0.8) {$ $};
\node (C1) at (0.3,1.2) {$ $};
\node (C2) at (1,1.2) {$ $};
\node (Comma) at (1.4,1.25) {.};
    \path (B1) edge [->]node [auto] {$\scriptstyle{s}$} (B2);
        \path (C1) edge [->]node [below] {$\scriptstyle{t}$} (C2);
  \end{tikzpicture}
\end{equation}
is oriented, i.e., the manifolds $U$ and $U\times_{\mathscr X} U$ and the maps $s$ and $t$ are oriented in a compatible way. In addition, we say that $\mathscr X$ is of {\em finite type} if $U$ and $U\times_{\mathscr X} U$ have finite good covers compatible with $s$ and $t$.
\end{definition}
Let $\mathscr X$ be a finite type smooth Deligne-Mumford differentiable stack. Then there is a well-defined theory of cohomology with compact supports for $\mathscr X$. In addition, there exists an \emph{integration} map
\begin{equation}\label{eq:integral}
\int_{\mathscr X}\colon H_{c}^{\dim(\mathscr X)}(\mathscr{X}) \to \mathbb{R}\ ,
\end{equation}
such that the induced pairing $H_{dR}^k(\mathscr X)\otimes H_{c}^{\dim(\mathscr X) - k}(\mathscr X)\to \mathbb{R}$ is perfect. If $\mathscr X$ is compact, $H_{dR}^k(\mathscr X)\simeq H_{c}^{k}(\mathscr X)$ for any $k$ and there is a well defined map \eqref{eq:integral} which induces a perfect pairing.

\subsection{Comparing algebraic, analytic, differentiable and topological stacks}

In \cite{art:lerman2010} a nice interpretation of differentiable stacks as Lie groupoids is thoroughly discussed, which  helped us to understand how the different types of stacks are related. 
 
Let $\mathscr{X}$ be a smooth separated Deligne-Mumford algebraic stack of finite type over $\mathbb{C}$ 
(cf.\ \cite[Sect.~4]{book:laumonmoret-bailly2000}). Then by \cite{art:keelmori1997}, there exists a {\em coarse moduli space} $\pi\colon \mathscr{X}\to X$. In general $X$ is a separated {\em algebraic space} of finite type over $\mathbb{C}$.
 
Let $\mathbf{AlgDM}$ be the 2-category of smooth separated Deligne-Mumford algebraic stacks of finite type over $\mathbb{C}$. 
As explained in \cite[Sect.~3.3]{art:behrendnoohi2006} (see also \cite{art:lerman2010}), 
it is equivalent to the weak 2-category  of groupoids up to Morita equivalence 

\begin{equation}
  \begin{tikzpicture}[xscale=2.8,yscale=-.7, ->, bend angle=25]
\node (A0_1) at (0,1) {$X_1$};
\node (A1_1) at (1.15,1) {$X_0$}; 
\node (B1) at (0.1,0.8) {$ $};
\node (B2) at (1,0.8) {$ $};
\node (C1) at (0.1,1.2) {$ $};
\node (C2) at (1,1.2) {$ $};
    \path (B1) edge [->]node [auto] {$\scriptstyle{s}$} (B2);
        \path (C1) edge [->]node [below] {$\scriptstyle{t}$} (C2);
  \end{tikzpicture}
\end{equation}
in the category of separated schemes of finite type over $\mathbb{C}$, where $X_0$ is smooth, $s, t$ are \'etale morphisms and 
\begin{equation}
(s,t)\colon X_1 \to X_0\times X_0
\end{equation}
is a proper unramified quasi-compact morphism (cf.\ \cite{book:laumonmoret-bailly2000}). Denote by $\mathbf{AnDM}$ (resp.\ $\mathbf{DiffDM}$, resp.\ $\mathbf{TopDM}$)
the 2-category of smooth Deligne-Mumford analytic stacks (resp.\ differentible stacks, resp.\ topological stacks).
The argument in \cite[Sect.~3.3]{art:behrendnoohi2006} proves that the 2-category of Deligne-Mumford topological stacks (resp.\ smooth Deligne-Mumford analytic stacks, resp.\ smooth Deligne-Mumford differentiable stacks)
is equivalent to the weak 2-category of groupoids up to Morita equivalence 
\begin{equation}
  \begin{tikzpicture}[xscale=2.8,yscale=-.7, ->, bend angle=25]
\node (A0_1) at (0,1) {$X_1$};
\node (A1_1) at (1.15,1) {$X_0$}; 
\node (B1) at (0.1,0.8) {$ $};
\node (B2) at (1,0.8) {$ $};
\node (C1) at (0.1,1.2) {$ $};
\node (C2) at (1,1.2) {$ $};
    \path (B1) edge [->]node [auto] {$\scriptstyle{s}$} (B2);
        \path (C1) edge [->]node [below] {$\scriptstyle{t}$} (C2);
  \end{tikzpicture}
\end{equation}
in the category of topological spaces (resp.\ complex manifolds, resp.\ differentiable manifolds), where the maps $s,t$ are local homeomorphisms (resp. local biholomorphisms, resp.\ local diffeomorphisms) and 
\begin{equation}
(s,t)\colon X_1 \to X_0\times X_0
\end{equation}
is a proper map with closed discrete fibers (resp.\ proper unramified finite map in the analytic and differential settings). Then one can define natural functors
\begin{equation*}
\mathbf{AlgDM} \xrightarrow{\mbox{-}^{\mathrm{an}}} \mathbf{AnDM}\xrightarrow{\mbox{-}^{\mathrm{diff}}} \mathbf{DiffDM} \xrightarrow{\mbox{-}^{\mathrm{top}}} \mathbf{TopDM}\ .
\end{equation*}
Moreover, these functors respect the coarse moduli space construction. We shall also denote by $\xrightarrow{\mbox{-}^{\mathrm{top}}}$, resp.\ $ \xrightarrow{\mbox{-}^{\mathrm{diff}}}$ any composition of these functors ending with $\xrightarrow{\mbox{-}^{\mathrm{top}}}$, resp.\ $\xrightarrow{\mbox{-}^{\mathrm{diff}}}$.
\begin{rem}
In the following, if $\mathscr{X}$ is a smooth Deligne-Mumford algebraic (analytic, differentiable) stack, its fundamental group is by definition the fundamental group of its underlying Deligne-Mumford topological stack. \hfill $\triangle$
\end{rem}

\bigskip\section{Differential geometry on smooth Deligne-Mumford analytic and differentiable stacks}\label{sec:diffstacks}

In this section we shall sketch the generalization of some basic results of differential geometry. 

\subsection{Vector bundles on smooth Deligne-Mumford differentiable stacks}

Let $\mathscr{X}$ be a smooth Deligne-Mumford differentiable stack.
\begin{definition}
A {\em $C^{\infty}$ complex vector bundle of rank $r$} on $\mathscr{X}$ is the following set of data:
\begin{itemize}
\item for any atlas $(U,u)$ a $C^{\infty}$ complex vector bundle $E_{U,u}$ of a rank $r$ on $U$,
\item for any morphism $(\varphi, \alpha)\colon (U,u)\to (V,v)$ of two atlases, where $\varphi\colon U\to V$ 
is a local diffeomorphism of differentiable manifolds and $\alpha\colon u\xrightarrow{\sim} v\circ \varphi$ is a 2-isomorphism, an isomorphism of $C^{\infty}$ complex vector bundles of rank $r$
\begin{equation*}
\vartheta_{\varphi,\alpha}\colon \varphi^\ast(E_{U,u})\xrightarrow{\sim} E_{V,v}
\end{equation*}
\end{itemize}
such that 
for any composition $(U,u)\xrightarrow{(\varphi,\alpha)} (V,v)\xrightarrow{(\psi, \beta)} (W,w)$, we have
\begin{equation*}
\varphi^\ast\vartheta_{\psi, \beta}\circ \vartheta_{\varphi,\alpha}=\vartheta_{\psi\circ \varphi, \varphi^\ast\beta\circ \alpha}\ .
\end{equation*}
\end{definition}
\begin{remark}
In \cite{art:behrendginotnoohixu2012} the authors give another definition of complex vector bundles on $\mathscr X$. Our definition is equivalent to their definition because 
of \cite[Prop.~3.2]{art:behrendginotnoohixu2012}.
\end{remark}

There are several operations on $C^{\infty}$ complex vector bundles on $\mathscr{X}$. We can indeed form the {\em trivial} vector bundle $V_\mathscr{X}:=V\times \mathscr{X}$ given any $\mathbb{C}$-vector space $V$, resp.\ the complex conjugate $\bar{E}$ or the dual ${E}^\vee$ of a vector bundle $E$, resp.\ the direct sum $E \oplus F$, the tensor product $E \otimes F$ or the bundle of morphisms $\mathrm{Hom}(E, F)$ of two vector bundles $E$ and $F$, resp.\ the pull back $f^\ast E$ of a complex vector bundle $E$ on $\mathscr{X}$ by a morphism $f\colon\mathscr{Y} \to \mathscr{X}$. Furthermore, one can define a morphism $\psi$ of $C^{\infty}$ complex vector bundles on $\mathscr{X}$ from  $E$ to $F$ as the data for any atlas $(U,u)$  of a morphism $\psi_{U,u}\colon E_{U,u}\to F_{U,u}$ commuting with $\vartheta_{\varphi,\alpha}$. This enables to define $C^{\infty}(\mathscr{X}, E)$, the {\em space of smooth sections of} $E$, as the set of morphisms from $\mathbb{C}_\mathscr{X}$ to $E$. 

A {\em connection} $\nabla$ on $E$ is the data of a connection $\nabla_{U,u}$ on $E_{U,u}$ for any atlas $(U,u)$ which is compatible to the  $\vartheta_{\varphi,\alpha}$.
In a similar vein, one can define riemannian metrics, hermitian metrics, Levi-Civita connection, 
principal bundles and connections, curvature tensors, \ldots.  The definitions are left to the reader, see however \cite{art:drydengordongreenwaldwebb2008} for orbifolds, 
and  \cite{art:behrendginotnoohixu2012}. The  theory of Chern-Weil forms of connections and the fact that, on smooth differentiable Deligne-Mumford stacks, 
they compute the rational Chern classes defined as in \cite{art:behrend2004}
is established in \cite{art:laurenttuxu:2007}. 

\subsection{Holonomy}\label{sec:holonomy}

Let $M$ be a differentiable manifold. Denote by $\mathrm{Conn}_M$ the groupoid whose objects are $C^{\infty}$
complex vector bundles over $M$ with a connection and whose arrows are {\em gauge equivalences} 
(vector bundle isomorphisms preserving the connection). Fix a $C^{\infty}$ complex vector bundle of rank $r$ with a connection $(E,\nabla)$
on a smooth connected Deligne-Mumford differentiable stack $\mathscr{X}$. Then, see \cite{art:laurenttuxu:2007},  we can define the holonomy functor
\begin{eqnarray*}
\mathrm{Hom}(M, \mathscr{X}) &\to& \mathrm{Conn}_M\ ,\\ 
f &\mapsto& (f^\ast E, f^\ast \nabla)\ .
\end{eqnarray*}
When $M=I$ is an open interval of the real numbers, $x_0,x_1\in I$ and $\delta\colon I\to \mathscr{X}$ is a {\em differentiable} path we can define the holonomy $h(\delta, x_0,x_1)$ along $\delta$ to be the natural linear transformation $\delta^\ast E_{x_0} \to \delta^\ast E_{x_1}$. If 
$\nabla$ is flat then $h(\delta, x_0,x_1)$ depends only on the differentiable homotopy class of $f\colon(I,\{x_0,x_1\})\to (\mathscr{X}, \{f(x_0), f(x_1)\})$. 

It is easy to define differentiable paths $[0,1] \to \mathscr{X}$ as {\em continuous} paths $[0,1] \to \mathscr{X}^{\mathrm{top}}$ extending differentiably to a slightly larger interval. Now if $\gamma\colon[0,1] \to \mathscr{X}^{\mathrm{top}}$ is a continuous path we can approximate it in the $C^0$-topology by differentiable paths. Since any homotopy of maps may be approximated by a differentiable homotopy in $C^1$-topology and homotopy classes are open in in $C^1$-topology, we conclude as in the classical case that the holonomy of the flat connection $\nabla$ on $E$ gives rise to a representation of the {\em Poincar\'e groupoid} of $\mathscr{X}$ hence to a group representation $\rho\colon\pi_1(\mathscr{X},x) \to GL(r,\mathbb{C})$. If the flat connection in question preserves a hermitian metric $h$, then $\rho$ takes values in the unitary group of $(E_x,h_x)$.  

Conversely, let $\tilde{\mathscr X}$ be a {\em universal cover}\footnote{The theory of covering stacks is developed in \cite[Sect.~18]{art:noohi2005}.} of
$\mathscr{X}$ and consider an atlas $(U,u)$ of $\mathscr X$. Then the fibred product $(\tilde{U}:=U\times_{\mathscr X} \tilde{\mathscr X}, \tilde{u})$ is 
an atlas of $\tilde{\mathscr X}$ and $\tilde{U}\to U$ is topological Galois covering space of  $U$ such that its deck transformations group is $\pi_1(\mathscr{X},y)$, 
where $y$ is point in $\mathscr{X}$.
Let $\rho\colon\pi_1(\mathscr{X},y)\to GL(r,\mathbb{C})$ be a linear representation.
Then $\pi_1(\mathscr{X},y)$ acts on $\tilde{U}\times \mathbb{C}^r$ by $\rho$ on the second factor and induces after passing to quotient
a $C^{\infty}$ complex vector bundle of rank $r$ with a flat connection $(E_{\rho_U}, \nabla_{\rho_U})$. 
This construction being functorial for morphisms of atlases, 
  one obtains a $C^{\infty}$ complex vector bundle on $\mathscr X$  of rank $r$ with a 
flat connection $(E_{\rho}, \nabla_{\rho})$ having $\rho$ as its holonomy. This construction is an equivalence of categories. 

\subsection{Metric geometry of Riemannian Deligne-Mumford stacks}

\begin{defin}
A {\em Riemannian Deligne-Mumford stack} is a pair  $(\mathscr{X}, g)$, where $\mathscr{X}$ is a smooth Deligne-Mumford differentiable stack and $g$ a Riemannian metric on it.\hfill $\oslash$
\end{defin}
Every smooth  Deligne-Mumford differentiable stack admits a Riemannian metric since we are assuming second countability. 
One can define geodesics on $\mathscr{X}$ or more generally harmonic mappings from manifolds to $(\mathscr{X}, g)$ in the usual
way  thanks to Proposition \ref{prop:localstructurestacks}. The infimal length of a path between two points 
of $\mathscr{X}^{\mathrm{top}}$ is the Riemannian distance $d$ between these two points and, since with our definitions the stacks are separated, this gives a distance function 
on $\mathscr{X}^{\mathrm{top}}$ inducing its topological structure (which is Hausdorff). The extension of the basic theory of geodesics 
is carefully carried out in \cite{art:guruprasadhaefliger2006}. 

For instance, the distance function permits to define the space of Lipschitz functions on $\mathscr{X}$ which can be characterized
by Rademacher's Theorem as almost everywhere differentiable functions with $g$-bounded differential. 
The space of continuous functions on $\mathscr{X}$ coincides with the space of continuous functions on $\mathscr{X}^{\mathrm{top}}$. 
However the space of Lipschitz (or of $C^{\infty}$) functions on $\mathscr X$ does not depend only on the structure of the underlying topological stack $\mathscr{X}^{\mathrm{top}}$, but it is a truly stacky invariant. 
For instance, if $G\subset GL(n, \mathbb{C})$ is a finite complex reflection group (e.g. $G=\{\pm 1 \} \in \mathbb{C}^\ast$) the distance function attached to a $G$-invariant hermitian form is just H\"older continuous with respect to a Riemannian metric in the usual sense on $\mathbb{C}^n /G \simeq \mathbb{C}^n$. In the example, $d_{[\mathbb{C}_w/\{\pm 1 \}]} (0,z)=|z|^{1/2}$ where $z=w^2$ is the natural coordinate on $\mathbb{C}_w=[\mathbb{C}_w/\{\pm 1 \}]^{\mathrm{top}}$. 

\subsection{Sheaves and vector bundles on smooth Deligne-Mumford analytic stacks}\label{sec:sheaves}
  
In this section we shall give definitions of (quasi)coherent analytic sheaves and vector bundles on a smooth Deligne-Mumford analytic stack 
which resembles the definition of (quasi)coherent sheaves given in \cite{art:vistoli1989} in the algebro-geometric setting. 

Let $\mathscr{X}$ be a smooth Deligne-Mumford analytic stack. 
\begin{definition}
A {\em quasicoherent sheaf} $\mathcal{F}$ on $\mathscr{X}$ is the following set of data:
\begin{itemize}
\item for any atlas $(U,u)$ a quasicoherent analytic sheaf $\mathcal{F}_{U,u}$ on $U$,
\item for any morphism $(\varphi, \alpha)\colon (U,u)\to (V,v)$ of two atlases, where $\varphi\colon U\to V$ is a local biholomorphism of complex manifolds and $\alpha\colon u\xrightarrow{\sim} v\circ \varphi$ is a 2-isomorphism, an isomorphism  
\begin{equation*}
\theta_{\varphi,\alpha}\colon \varphi^\ast\mathcal{F}_{U,u}\xrightarrow{\sim} \mathcal{F}_{V,v}
\end{equation*}
\end{itemize}
such that 
for any composition $(U,u)\xrightarrow{(\varphi,\alpha)} (V,v)\xrightarrow{(\psi, \beta)} (W,w)$, we have
\begin{equation*}
\varphi^\ast\theta_{\psi, \beta}\circ \theta_{\varphi,\alpha}=\theta_{\psi\circ \varphi, \varphi^\ast\beta\circ \alpha}\ .
\end{equation*}
A {\em coherent sheaf} is a quasicoherent sheaf $\mathcal{F}$ such that all $\mathcal{F}_{U,u}$ are coherent. A {\em locally free sheaf} is a coherent sheaf $\mathcal{F}$ such that all $\mathcal{F}_{U,u}$ are locally free. A morphism $f\colon \mathcal{F}\to \mathcal{F}'$ is a collection of morphisms $f_{U,u}\colon \mathcal{F}_{U,u}\to \mathcal{F}_{U,u}'$ such that for any morphism $(\varphi, \alpha)\colon (U,u)\to (V,v)$ of two atlases we have $\varphi^\ast(f_{V,v})\circ \theta_{\varphi,\alpha} = \theta_{\varphi,\alpha}'\circ f_{U,u}$.
\end{definition}  
\begin{example}
The structure sheaf $\mathcal{O}_{\mathscr{X}}$ of $\mathscr{X}$ is defined by $(\mathcal{O}_{\mathscr{X}})_{U,u}=\mathcal{O}_U$ for any atlas $(U,u)$, with the obvious isomorphisms. The sheaf of differential $\Omega^1_{\mathscr{X}}$ of $\mathscr{X}$ is defined by $(\Omega^1_{\mathscr{X}})_{U,u}=\Omega^1_U$. Since for any morphism $(\varphi, \alpha)\colon (U,u)\to (V,v)$ of two atlases, the morphism $\varphi$ is a local biholomorphism, there is a natural isomorphism $\Omega^1_U\simeq \varphi^\ast(\Omega^1_V)$.
\end{example}
A {\em holomorphic vector bundle of rank $r$} on $\mathscr{X}$ is a $C^{\infty}$-complex vector bundle $E$ over $\mathscr{X}$ such that all $E_{U,u}$ are holomorphic vector bundles and all $\vartheta_{\varphi,\alpha}$ are isomorphisms of holomorphic vector bundles. As before in the differential setting, also in this setting one can perform the usual operations on vector bundles. 

Let $E$ be a rank $r$ holomorphic vector bundle on $\mathscr{X}$. For any atlas $(U,u)$ we denote by $\mathcal{E}_{U,u}$ the sheaf of holomorphic sections of the vector bundle $E_{U,u}$. For any morphism $(\varphi, \alpha)\colon (U,u)\to (V,v)$ of two atlases, the isomorphism $\vartheta_{\varphi,\alpha}$ induces an isomorphism 
\begin{equation}
\theta_{\varphi,\alpha}\colon \varphi^\ast\mathcal{E}_{U,u}\xrightarrow{\sim} \mathcal{E}_{V,v}\ .
\end{equation}
In this way, we define a coherent sheaf  $\mathcal{E}$ on $\mathscr{X}$, the {\em sheaf of sections of the vector bundle} $E$.
\begin{remark}
The functor $-^{\mathrm{an}}$ sends coherent sheaves to coherent analytic sheaves. The functor $-^{\mathrm{diff}}$ sends holomorphic vector bundles to $C^\infty$ complex vector bundles.
\end{remark}

Let $\mathscr{X}$ be a smooth Deligne-Mumford analytic stack. Then it is easy to see that there is a one-to-one correspondence between holomorphic vector bundles of rank $r$ on $\mathscr{X}$ and locally free sheaves of rank $r$ on it. Assume moreover that $\mathscr{X}$ is of the form $[Z/G]$, where $Z$ is a complex manifold and $G$ a complex Lie group acting on it. By the same argument as in \cite[Example~12.4.6]{book:laumonmoret-bailly2000}, the category of coherent sheaves on $\mathscr{X}$ is equivalent to the category of $G$-equivariant coherent sheaves on $Z$. Similarly, the category of locally free sheaves on $\mathscr{X}$ is equivalent to the category of $G$-equivariant locally free sheaves on $Z$ (this follows from the descent result for locally free sheaves with respect to fppf morphisms, \cite[Prop.~A.11]{book:algebraicstacks}). Since the category of $G$-equivariant locally free sheaves on $Z$ is equivalent to the category of $G$-equivariant holomorphic vector bundles on $Z$ (see, e.g., \cite[
Sect.~5.1]{book:chrissginzburg2010}), the category of holomorphic vector bundles of rank $r$ on $\mathscr{X}$ is equivalent to the category of $G$-equivariant holomorphic vector bundles of rank $r$ on $Z$.

In the following we shall call {\em hermitian bundle} a holomorphic vector bundle on $\mathscr{X}$ with a hermitian metric on it. Then if $\mathscr{X}$ is a smooth Deligne-Mumford analytic stack of the form $[Z/G]$, the category of hermitian bundles on $\mathscr{X}$ is equivalent to the category of $G$-equivariant hermitian vector bundles on $Z$. Thus if $E$ is a holomorphic vector bundle on a smooth Deligne-Mumford analytic stack, by Proposition \ref{prop:localstructurestacks} and a partition of unity argument on $\mathscr{X}^{\mathrm{top}}$, one can always construct a hermitian metric on $E$.

The usual constructions of complex differential geometry - see \cite[Chap. V]{book:demailly} - carry over to smooth Deligne-Mumford analytic stacks\footnote{One can mimic the arguments in \cite{art:behrend2004}: one can define everything and the level of groupoid presentations of the stack and then prove that it is invariant under Morita equivalence.}. A smooth Deligne-Mumford analytic stack carries a complexified tangent bundle $T_\mathscr{X}^{\mathbb{C}}= T^{1,0}_\mathscr{X} \oplus T^{0,1}_\mathscr{X}$ from which we may form the usual vector bundles $\Omega^k_\mathscr{X}$ of $k$-forms and $\Omega^{p,q}_\mathscr{X}$ of $(p,q)$-forms. A $(1,1)$ form $\omega$ on $\mathscr{X}$ is {\em positive} if so are its representatives $\omega_{U,u}$ for any atlas $(U,u)$ of $\mathscr{X}$. A hermitian vector bundle $(E,h)$ carries a canonical Chern connection $\nabla^h$ with a curvature form $i\Theta(E,h) \in C^{\infty}(\mathscr{X}, \mathrm{Herm}(E,h)\otimes \Omega^{1,1}_\mathscr{X})$. The curvature of a hermitian line 
bundle $(L,h)$ is just a real $(1,1)$ form on $\mathscr{X}$.

\subsection{K\"ahler Deligne-Mumford analytic stacks}

At a slightly deeper level, the De Rham, Dolbeault \cite[Chap.~IV, Sect.~6]{book:demailly} and Hodge \cite[Chap.~VI, Thm.~3.16 and 3.17]{book:demailly} 
isomorphism theorems are valid for smooth Deligne-Mumford analytic stacks. If we follow the proof of these basic results \cite{book:demailly}, we see that 
the only non obvious points are integration of top dimensional compactly supported forms on an oriented smooth Deligne-Mumford differentiable stack 
(so that we can define the $L_2$ norm on the space of sections of a metrized vector bundle on an oriented\footnote{One needs to use densities for non-orientable stacks.} Riemannian Deligne-Mumford differentiable stack), 
Stokes' theorem (so that we can perform integration by parts) and acyclicity of sheaves of $(p,q)$-forms (cf.\ Section \ref{sec:analdiffstacks} 
for the first two results, the latter descends from the analogous statement in the case of complex manifold simply by using a groupoid presentation of the stack). 

The literature on the Grothendieck-Riemann-Roch theorem is not as complete as one would like, even for complex manifolds, see however \cite{art:toen1999} for proper representable morphisms of smooth Deligne Mumford algebraic stacks over the complex numbers. Fortunately, this is not used here. 

Say a smooth Deligne-Mumford analytic stack is {\em K\"ahler} if it carries a closed positive $(1,1)$-form. This is equivalent to require that the coarse moduli space is K\"ahler. The Lefschetz package, the $\partial\bar \partial$-lemma hold on compact K\"ahler Deligne-Mumford analytic stacks and so does the Kodaira-Akizuki-Nakano theorem. An application is an analytic version of  Olsson-Matsuki's proof of the Kawamata-Viehweg vanishing theorem in characteristic zero \cite{art:matsukiolsson2005}. Yau's solution of the Calabi conjecture \cite{art:yau1978} is also valid on stacks. For orbifolds, an early reference is \cite{art:kobayashi1985}, see also \cite{art:eyssidieuxguedjzeriahi2009}. Although we will not use this, we note that the characterization of the K\"ahler cone, extends to compact K\"ahler Deligne-Mumford analytic stacks. 

\begin{theorem}\cite{art:demaillypaun2004}
Let $\mathscr{X}$ be a $n$-dimensional compact K\"ahler Deligne-Mumford analytic stack and let $\omega$ be K\"ahler form on $\mathscr{X}$. A class $\{\alpha\}\in H^{1,1}(\mathscr{X}, \mathbb{R})$ is represented by a K\"ahler form iff for every closed irreducible analytic substack $\mathscr{Z} \subset \mathscr{X}$ one has $\{\alpha\}^{\dim (\mathscr{Z})}.[\mathscr{Z}]>0$ and $\displaystyle\int_{\mathscr{X}}\{\alpha\}^{n-p}\,\{\omega\}^p >0$ for all $1\le p \le n-1$. 
\end{theorem}
We will not need the full strength of this statement and will only sketch the proof of the consequence we shall use. 

\begin{definition} 
Let $\mathcal{L}$ be a holomorphic line bundle on a smooth Deligne-Mumford analytic stack $\mathscr{X}$. Let $\mathscr{Z} \subseteq \mathscr{X}$ be a closed substack of $\mathscr{X}$. We say $\mathcal{L}$ is {\em positive} on $\mathscr{Z}$ if there is an open analytic substack $\mathscr{U}$ such that $\mathscr{Z} \subset \mathscr{U}$ and a hermitian metric $h$ on $\mathcal{L}$ such that $i\Theta(L_{\vert\mathscr{U}},h_{\vert\mathscr{U}})$ is positive. 
\end{definition}

\begin{proposition}\label{prop:positiveness}
Let $\mathcal{L}$ be a holomorphic line bundle on a smooth Deligne-Mumford analytic stack $\mathscr{X}$. Assume there exists $n\in \mathbb{N}\setminus\{0\}$ such that  
$\mathcal{L}^{\otimes n}\simeq \pi^\ast(M)$ where $\pi\colon \mathscr{X}\to X$ is the natural map to the coarse moduli space $X$ of $\mathscr{X}$ and $M$ 
is a holomorphic line bundle on $X$. Let $Z\subseteq X$ be a compact analytic subspace. Then $\mathcal{L}$ is positive on $\mathscr{Z}:=Z\times_X \mathscr{X}$
 if and only if $M$ is ample on $Z$ (which is then a projective algebraic variety). 
\end{proposition}
\proof 
The corresponding statement for analytic spaces was proved in P\u aun's thesis, see \cite[Prop.~3.3]{art:demaillypaun2004} for the required level of generality. By Proposition \ref{prop:localstructurestacks}, it rests on gluing and regularization techniques for quasi-plurisubharmonic functions \cite{art:demailly1990} that are local in nature and compatible with the action of a finite group. 
\endproof
\begin{remark}
The previous proposition holds also for Deligne-Mumford analytic stacks that need not be smooth.
\end{remark}
\begin{remark}
A positive line bundle in the above sense needs not be ample in the sense of \cite{art:rossthomas2011} because this latter condition requires that the stack has cyclic isotropy groups and the isotropy representations of the  line bundle are faithful. See \cite{art:kresch2009} for related issues. 
\end{remark}

\section{Hermite-Einstein metrics on compact K\"ahler Deligne-Mumford analytic stacks}\label{sec:hermite-einstein}

\subsection{Uhlenbeck-Yau theorem}
Let $\mathscr{X}$ be a $n$-dimensional connected compact K\"ahler Deligne-Mum\-ford analytic stack and let $\omega$ be K\"ahler form on $\mathscr{X}$. 

Let $\mathcal{E}$ be a coherent sheaf on $\mathscr{X}$. We call {\em rank} of $\mathcal{E}$ the zero degree part of the Chern character of its K-theory class $[\mathcal{E}]$.
\begin{definition}
A holomorphic vector bundle $\mathcal{E}$ on $\mathscr{X}$ is {\em $\omega$-stable} if every coherent subsheaf $\mathcal{F} \subset \mathcal{E}$ satisfies 
 $\mu(\mathcal{F})<\mu(\mathcal{E})$, where $\mu(\mathcal{F})$ is the {\em $\omega$-slope} of $\mathcal{F}$ defined as
\begin{equation*}
\mu(\mathcal{F}):=\frac{\displaystyle\int_{\mathscr{X}} \mathrm{c}_1(\mathcal{F}).\{\omega\}^{\dim(\mathscr{X})-1}}{\mathrm{rk}(\mathcal{F})}\ .
\end{equation*} 
\end{definition}
The proof of the celebrated Uhlenbeck-Yau theorem carries over to K\"ahler Deligne-Mumford analytic stacks. 

\begin{theorem} \label{thm:uy} \cite{art:uhlenbeckyau1986}
Let $\mathscr{X}$ be a $n$-dimensional connected compact K\"ahler Deligne-Mumford analytic stack and let $\omega$ 
be K\"ahler form on $\mathscr{X}$. Let $\mathcal E$ be a $\omega$-stable holomorphic vector bundle on $\mathscr{X}$. Then it carries a smooth Hermite-Einstein metric $h$ namely a hermitian metric satisfying the equation: 
\begin{equation*}
\Lambda i\Theta(\mathcal E,h)= 2\pi \mu(\mathcal{E}) \mathrm{Id}_{\mathcal E}\ .
\end{equation*}
\end{theorem}
 
\subsection{Comparison with Nironi's stability condition}

Let $(\mathscr{P},\mathcal{O}_{\mathscr{P}}(1))$ be a $n$-dimensional polarized projective cyclic orbifold in the sense of \cite{art:rossthomas2011}. 
One may take 
$\mathscr{P}:=\mathscr{P}(a_0, \ldots, a_n)$ the $n$-dimensional weighted projective stack for some coprime positive integers $a_0, \ldots, a_n$ with $n\geq 1$. 
Then for some $m$, $\mathcal{G}:=\oplus_{i=1}^m \mathcal{O}_{\mathscr{P}}(i)$ is a generating sheaf \cite[Sect.~5.2]{art:kresch2009}. 

In \cite{art:rossthomas2011} the authors define a Kahler form $\omega$ induced by the curvature of a metric on $\mathcal{O}_{\mathscr{P}}(1)$. Then for any $n$-dimensional coherent sheaf $\mathcal{E}$ on $\mathscr{X}$ we have
\begin{equation}
\frac{1}{\mathrm{rk}(\mathcal{E})}\int_{\mathscr{P}}\mathrm{c}_1(\mathcal{E})\cdot \{\omega\}^{n-1}=\frac{1}{\mathrm{rk}(\mathcal{E})}\int_{\mathscr{P}}\mathrm{c}_1(\mathcal{E})\cdot \mathrm{c}_1(\mathcal{O}_{\mathscr{P}}(1))^{n-1}\ .
\end{equation}
By using \cite[Prop.~3.18]{art:nironi2008} one can prove
\begin{equation}
\frac{1}{\mathrm{rk}(\mathcal{E})}\int_{\mathscr{P}}\mathrm{c}_1(\mathcal{E})\cdot \mathrm{c}_1(\mathcal{O}_{\mathscr{P}}(1))^{n-1}=\frac{1}{m^n}\mu_{\mathcal{G}}(\mathcal{E})\ ,
\end{equation}
where $\mu_{\mathcal{G}}(\mathcal{E})$ is the $\mathcal{G}$-slope introduced by Nironi. Therefore the $\omega$-stability condition is equivalent to the $\mu_{\mathcal{G}}$-stability condition of Nironi.

\section{Hermite-Einstein metrics on some noncompact K\"ahler Deligne-Mumford analytic stacks}\label{sec:hermite-einstein2}
 
In this section we shall prove a variant of \cite{art:bando1993} for stacks.

\subsection{Deleted neighborhoods of smooth divisors}
Let  $(\mathscr{X},g)$ be a  Riemannian Deligne-Mumford stack. 
Given an atlas $(U,u)$ of $\mathscr{X}$ one can construct a groupoid action
of  
\begin{equation}
  \begin{tikzpicture}[xscale=2.8,yscale=-.7, ->, bend angle=25]
\node (A0_1) at (-0.2 ,1) {$U\times_{\mathscr{X}}U$};
\node (A1_1) at (1.15,1) {$U$}; 
\node (B1) at (0.1,0.8) {$ $};
\node (B2) at (1,0.8) {$ $};
\node (C1) at (0.1,1.2) {$ $};
\node (C2) at (1,1.2) {$ $};
\node (Comma) at (1.4,1.25) {};
    \path (B1) edge [->]node [auto] {$\scriptstyle{s}$} (B2);
        \path (C1) edge [->]node [below] {$\scriptstyle{t}$} (C2);
  \end{tikzpicture}
\end{equation}
on the 
tangent space to $U$  and the resulting quotient stack $\pi\colon \mathrm{Tot} (T \mathscr{X})\to \mathscr{X}$ is a differentiable Deligne-Mumford stack over $\mathscr{X}$ which is the total space of the tangent vector bundle $T \mathscr{X}$ to $\mathscr{X}$. One can generalize this construction to any (real) vector bundle $E$ on $\mathscr X$.

The zero section of any (real) vector bundle $E$ defines a substack  $Z(E) \to  \mathrm{Tot}(E)$ of its total space. 
Moreover, $\mathrm{Tot}(E)\setminus Z(E)$ is representable (i.e.: it is a manifold) if and only if the  isotropy action
at all points of $\mathscr{X}$ has no non zero fixed vectors. 
Given a metric on $E$ one can construct similarily the unit sphere bundle $S(E)$ as a smooth substack of $\mathrm{Tot}(E)\setminus Z(E)$ which is a locally trivial sphere fibration
over $\mathscr{X}$. Moreover, the long exact sequence of homotopy groups holds true \cite{art:noohi2005}. 
In particular, $\mathrm{Tot}(T \mathscr{X})\setminus \mathscr{X}$ is an honest manifold (i.e.: it is representable) if and only if the tangent isotropy actions have no non zero fixed vectors. 
One also has an equivalence $\mathrm{Tot}(E)\setminus Z(E) \to S(E) \times \mathbb{R}_{>0}$ the second factor $\rho\colon \mathrm{Tot}(E)\to \mathbb{R}_{>0}$ being the norm function. 

Given $i\colon \mathscr{S} \to \mathscr{X}$ 
a closed smooth substack,
the normal bundle $\mathscr{N}_{\mathscr{S}|\mathscr{X}}$ over $\mathscr{S}$ can be constructed as $i^* T\mathscr{X} / T \mathscr{S}$. 
Then, using the exponential map in normal directions, one constructs an open differentiable substack $\mathscr{V} \subset \mathscr{X}$ containing $\mathscr{S}$ as a closed substack
such that the pair $(\mathscr{V},\mathscr{S})$  is equivalent to $(\mathrm{Tot}(\mathscr{N}_{\mathscr{S}|\mathscr{X}}), \mathscr{S})$. 

If $\mathscr{X}$ is now a smooth Deligne-Mumford analytic stack and $\mathscr{D}\subset \mathscr{X}$ be a smooth integral codimension one closed substack. Then $\mathscr{N}_{\mathscr{D}^{\mathrm{diff}}|\mathscr{X}^{\mathrm{diff}}}$ is the smooth bundle underlying $\mathcal{O}_{\mathscr{X}}(\mathscr{D})$ and there is an open differentiable substack $\mathscr{V} \subset \mathscr{X}^{\mathrm{diff}}$ containing $\mathscr{D}^{\mathrm{diff}}$ as a retract 
and a homotopy equivalence $\mathscr{U}_{\mathscr{D}}:=S(\mathscr{N}_{\mathscr{D}^{\mathrm{diff}}|\mathscr{X}^{\mathrm{diff}}}) \to \mathscr{V}\setminus \mathscr{D}^{\mathrm{diff}}$ such that
the resulting map $\mathscr{U}_{\mathscr{D}}\to \mathscr{D}^{\mathrm{diff}}$ is a locally trivial $S^1$-fibration. In particular one has an exact sequence:
\begin{equation}\label{eq:shortexact}
1 \to C_{\mathscr{D}} \to \pi_1(\mathscr{U}_{\mathscr{D}}) \to \pi_1(\mathscr{D}) \to 1\ ,
\end{equation}
where $C_{\mathscr{D}}$ is a normal cyclic subgroup whose order may take any non negative integer value or $+\infty$. If $\mathcal{O}_{\mathscr{X}}(\mathscr{D})$ is ample on $\mathscr{D}$, then $\mathscr{U}_{\mathscr{D}}$ is a complex manifold.

\subsection{Asymptotically flat metrics on deleted neighborhoods}

Assume  now  $(\mathscr{X}, \omega_0)$ is K\"ahler Deli\-gne-Mumford stack and $\mathscr{D}\subset \mathscr{X}$ is a smooth divisor. 
A riemannian metric $g_c$ on $\mathscr{X}\setminus{\mathscr{D}}$ has cone-like singularities if there exists a metric $ds^2$
on $\mathscr{U}_{\mathscr{D}}$ such that the restriction of $g_c$ to $\mathscr{V}\setminus \mathscr{D}^{\mathrm{diff}}$ is asymptotic
to $dr^2+r^2ds^2$ in the sense made precise by \cite[Sect.~1, Def.~1]{art:bando1993} on manifolds 
where $r$ is the distance function to any given point (note that $r\sim \rho^{-1}$). The curvature of such metric $g_c$ decays to zero 
when approaching $\mathscr{D}$. 

If  $\mathcal{O}_{\mathscr{X}}(\mathscr{D})|_{\mathscr{D}}$ is positive, there exists an hermitian
metric $h_{\mathscr{D}}$ on $\mathcal{O}_{\mathscr{X}}(\mathscr{D})$ such that its curvature form is positive definite on a neighbourhood of $\mathscr{D}$. 
As in \cite[Sect.~1]{art:bando1993}, we can construct a cone-like complete metric $\omega_c$ on $\mathscr{X}\setminus \mathscr{D}$ by the following formula:
\begin{equation}
 \omega_c:=\sqrt{-1}\partial\bar\partial \frac{1}{a}\exp(a\log\|\sigma_\mathscr{D} \|^{-2}_{h_{\mathscr{D}}}) + C \omega_0\ ,
\end{equation}
where $\sigma_{\mathscr{D}} \in H^0(\mathscr{X}, O_{\mathscr{X}}(\mathscr{D}))$ is the {\em tautological section} of $\mathscr{D}$, $a$ is an arbitrary positive number
and $C$ is a big enough positive real number.

\begin{definition}
 A holomorphic hermitian vector bundle $(\mathcal E,h)$ on $\mathscr{X}\setminus \mathscr{D}$ is asymptotically flat
if $|\Theta(\mathcal E_{\vert(X^{\circ},\omega_c)},h)|_{g_c}=O(r^{-2-\epsilon})$ as $r\to \infty$
for some $\epsilon >0$.
\end{definition}

Restricting the Chern connection of $(\mathcal E,h)$ to $\{r=R\}$ 
for $R\gg 1$ we get a hermitian connection $D_R$ on the pull-back of $\mathcal E$ to $\mathscr{U}_{\mathscr{D}}$. 

The asymptotic flatness condition gives immediately (cf.\ \cite[Sect.~1, Def.~3]{art:bando1993}):

\begin{lemma}
The holonomy of $D_R$ along any path has  a limit when $R\to \infty$ which is homotopy invariant.   
\end{lemma}

This defines a unitary representation of $\pi_1(\mathscr{U}_{\mathscr{D}})$ and 
we say $(\mathcal E,h)$ has {\em trivial holonomy at infinity} if the representation is trivial when restricted to the group $C_{\mathscr{D}}$ defined by \eqref{eq:shortexact}. 

\subsection{Bando's Instanton Theorem}

We are now ready to state an analog of \cite[Thm.~1]{art:bando1993}:

\begin{theorem}\label{thm:bando}
Let $(\mathscr{X}, \omega_0)$ be a $n$-dimensional ($n\geq 2$) connected compact K\"ahler Deligne-Mumford stack and let $\mathscr{D}\subset \mathscr{X}$ be a smooth  divisor such that $\mathcal{O}_{\mathscr{X}}(\mathscr{D})$ is ample and positive on $\mathscr{D}$. Let $\omega_c$ be a K\"ahler metric  with cone-like singularities on $X^{\circ}:=\mathscr{X}\setminus \mathscr{D} $. Then a holomorphic vector bundle $\mathcal E$ on $\mathscr{X}$ is such that $\mathcal E_{\vert\mathscr{D}}$ can be endowed with a flat unitary connection $\nabla$ iff there exists an Hermite-Einstein asymptotically flat vector bundle $(\mathcal{E}',h')$ on $(X^{\circ},\omega_c)$ having trivial holonomy at infinity which induces $\nabla$ and $\mathcal{E}$ is an extension of $\mathcal{E}'$.
\end{theorem}
\proof
 The previous section gives one direction of the first statement. 
Bando's proof of the much harder converse  applies mutatis mutandis: extend the flat hermitian metric to a deleted neighborhood of $\mathscr{D}$ by using a retraction, then 
extend it to $\mathscr{X}$
and flow it to an
asymptotically flat Hermite-Einstein metric. The extension theorem, being local, also follows in the same way. 
\endproof

\begin{corollary}[{\cite[Cor.~2]{art:bando1993}}]\label{cor:bando}
If furthermore $X^{\circ}$ is an ALE {\em manifold} then
we can replace the curvature decay condition by $\displaystyle\int_\mathscr{X} |\Theta(\mathcal{E}',h')|^n <+\infty$; it is equivalent to $|\Theta(\mathcal{E}',h')|=O(r^{-(2n-\epsilon)})$ for any $\epsilon >0$. 
\end{corollary}

\bigskip\section{An application: orbifold compactifications of ALE spaces of type $A_{k-1}$}\label{sec:application}

In this section we describe an explicit application of Theorem \ref{thm:bando}.

\subsection{Homotopy theory of spherical curves}

In this section we shall recall some homotopical properties of spherical Deligne-Mumford curves from \cite{art:behrendnoohi2006}.

For any pair of integers $m,n\geq 1$, we call {\em weighted projective line of type $(m,n)$} the smooth Deligne-Mumford analytic stack $\mathscr{P}(m,n):=[\mathbb{C}^2\setminus \{0\}/\mathbb{C}^\ast]$, where the action of $\mathbb{C}^\ast$ on $\mathbb{C}^2\setminus \{0\}$ is given by $t\cdot (x,y):=(t^mx,t^ny)$ for any $t\in\mathbb{C}^\ast$ and $(x,y)\in\mathbb{C}^2\setminus \{0\}$. Note that $\mathscr{P}(1,1)\simeq \mathbb{P}^1$. The stack $\mathscr{P}(m,n)$ has at most two orbifold points and its coarse moduli space is $\mathbb{P}^1$. Obviously, $\mathscr{P}(m,n)\simeq \mathscr{P}(n,m)$. Moreover, a weighted projective line is an orbifold if an only if $m$ and $n$ are relatively coprime. We call these {\em orbifold weighted projective lines}. The weighted projective line $\mathscr{P}(m,n)$ is a $\mu_d$-gerbe over $\mathscr{P}(\frac{m}{d},\frac{n}{d})$, where $d=\mathrm{gcd}(m,n)$. As explained in \cite[Sect.~4.3]{art:behrendnoohi2006}, the fundamental group $\pi_{1}(\mathscr{P}(m,n))$ is trivial.

For any pair of integers $m,n\geq 1$, a {\em football of type $(m,n)$} is an orbifold $\mathscr{F}(m,n)$ whose coarse moduli space is $\mathbb{P}^1$ and has two orbifold points of order $m$ and $n$ at $0$ and $\infty$, respectively. When $m$ and $n$ are relatively prime, $\mathscr{F}(m,n)$ is naturally isomorphic to $\mathscr{P}(m,n)$. 

\begin{definition}
Let $\mathscr{D}$ be a one-dimensional smooth Deligne-Mumford analytic stack. We say that $\mathscr{D}$ is a spherical Deligne-Mumford curve if its universal cover is $\mathscr{P}(m,n)$ for some positive integer numbers $m,n$.
\end{definition}
\begin{remark}
Let $\mathscr{D}$ be one-dimensional smooth Deligne-Mumford analytic stack. By \cite[Prop.~4.6]{art:behrendnoohi2006} there exists a one-dimensional analytic orbifold $\tilde{\mathscr D}$ and a finite group $H$ such that $\mathscr{D}$ is a $H$-gerbe over $\tilde{\mathscr D}$. By \cite[Prop.~7.4]{art:behrendnoohi2006} $\mathscr{D}$ is spherical if and only if $\tilde{\mathscr{D}}$ is spherical.
\end{remark}
\begin{proposition}{\normalfont \cite[Prop.~5.5]{art:behrendnoohi2006}.}\label{prop:fundamentalgroup-orbifoldcurves}
Let $\mathscr{D}$ be a spherical Deligne-Mumford curve which is an orbifold with two orbifold points. Then $\mathscr{D}$ is isomorphic to a football $\mathscr{F}(m,n)$ of type $(m,n)$ for some positive integers $m,n$. The fundamental group of $\mathscr{D}$ is $\mathbb{Z}_d$, where $d=\mathrm{gcd}(m,n)$ and its universal cover is $\mathscr{P}(\frac{m}{d},\frac{n}{d})$ on which $\mathbb{Z}_d$ acts by rotations.
\end{proposition} 
\begin{proposition}\label{prop:fundamentalgroup-analyticcurves}
Let $\mathscr{D}$ be a one-dimensional smooth Deligne-Mumford analytic stack, which is an $H$-gerbe over a football $\mathscr{F}(m,n)$. Then we have the exact sequence
\begin{equation}
1\to \mathbb{Z}_d\to H\to \pi_1(\mathscr{D})\to \pi_1(\mathscr{F}(m,n))\to 1\ ,
\end{equation}
where $d=\mathrm{gcd}(m,n)$. Moreover, the universal cover of $\mathscr{D}$ is isomorphic to $\mathscr{P}(m,n)$.
\end{proposition}
\proof
By \cite[Prop.~7.6-(i)]{art:behrendnoohi2006}, the fundamental group of $\mathscr{D}$ fits into the exact sequence
\begin{equation*}
1\to \ker(H\to \pi_1(\mathscr{D}))\to H\to \pi_1(\mathscr{D})\to \pi_1(\mathscr{F}(m,n))\to 1\ .
\end{equation*}
On the other hand, as explained in \cite[Sect.~9]{art:behrendnoohi2006}, $\ker(H\to \pi_1(\mathscr{D}))\simeq \pi_2(\mathrm{PGL}(m,n))$ (for a definition of $\mathrm{PGL}(m,n)$, see \cite[Sect.~8]{art:behrendnoohi2006}). Finally, by \cite[Prop.~8.5 and 8.6]{art:behrendnoohi2006}, we get
\begin{equation*}
\pi_2(\mathrm{PGL}(m,n))=\left\{
\begin{array}{ll}
\mathbb{Z}_d & \mbox{if $m\neq n$ and $d=\mathrm{gcd}(m,n)$}\ ,\\
\mathbb{Z}_m & \mbox{if $m=n$}\ .
\end{array}
\right.
\end{equation*}
The second statement follows from \cite[Prop.~7.6-(ii)]{art:behrendnoohi2006}.
\endproof

\subsubsection{Spherical Deligne-Mumford curves with universal cover $\mathscr{P}(d,d)$}\label{sec:fundamental-cover}

Let $\mathscr{D}$ be a one-dimensional smooth Deligne-Mumford analytic stack with universal cover $\mathscr{P}(d,d)$. By \cite[Cor.~9.8]{art:behrendnoohi2006}, $\mathscr{D}$ is a global quotient stack of the form $[\mathbb{C}^2\setminus \{0\}/E]$ where $E$ is a central extension of a discrete group $\Gamma$ by $\mathbb{C}^\ast$ and the action of $E$ on $\mathbb{C}^2\setminus \{0\}$ is given by a $\mathbb{C}$-representation $\rho\colon E\to \mathrm{GL}(2,\mathbb{C})$ of $E$, which fits into the commutative triangle
\begin{equation}\label{eq:triangle}
  \begin{tikzpicture}[xscale=1.5,yscale=-1.2]
    \node (A0_1) at (1, 0) {$\mathbb{C}^\ast$};
    \node (A1_0) at (0, 1) {$E$};
    \node (A1_2) at (2, 1) {$\mathrm{GL}(2,\mathbb{C})$};
    \path (A0_1) edge [->]node [auto] {$\scriptstyle{}$} (A1_0);
    \path (A0_1) edge [->]node [auto] {$\scriptstyle{(t^d,t^d)}$} (A1_2);
    \path (A1_0) edge [->]node [auto] {$\scriptstyle{\rho}$} (A1_2);
  \end{tikzpicture} 
\end{equation}  
As explained in \cite[Sect.~9.2]{art:behrendnoohi2006}, the fundamental group of $\mathscr{D}$ is $\Gamma$. 
\begin{remark}\label{rem:flatconnection-sphericalcurves}
The Picard group $\operatorname{Pic}(\mathscr{D})$ of $\mathscr{D}$ consists of $E$-equivariant line bundles on $\mathbb{C}^2\setminus\{0\}$, i.e.,
characters $E\to \mathbb{C}^\ast$ (cf.\ Section \ref{sec:sheaves}). Let $\rho\colon \Gamma\to \mathbb{C}^\ast$ be a one-dimensional representation of $\Gamma\simeq \pi_1(\mathscr{D})$. 
Then the composition of morphisms $E\to \Gamma\xrightarrow{\rho} \mathbb{C}^\ast$ defines a line bundle $L_\rho$ on $\mathscr{D}$. 
Moreover, $L_\rho$ is endowed by an hermitian metric and a unitary flat connection associated with $\rho$ (cf.\ Section \ref{sec:holonomy}). 
\end{remark}

\subsection{Orbifold compactification}

In this section we describe a construction in \cite{art:bruzzopedrinisalaszabo2013} of a compactification of the minimal resolution of the $A_{k-1}$ toric singularity of $\mathbb{C}^2/\mathbb{Z}_k$, which turns out to be a projective toric orbifold.

\subsubsection*{Normal compactification}\label{sec:minimalresolution}

Let $k\geq 2$ be an integer and denote by $\mu_k$ the group of $k$-th roots of unity in $\mathbb{C}$. A choice of a primitive $k$-th root of unity $\omega$ defines an isomorphism of groups $\mu_k\simeq \mathbb{Z}_k$. We define an action of $\mu_k\simeq \mathbb{Z}_k$ on $\mathbb{C}^2$ as
\begin{equation}
\omega\cdot (z_1,z_2):= (\omega z_1, \omega^{-1} z_2)\ .
\end{equation}
The quotient $\mathbb{C}^2/\mathbb{Z}_k$ is a normal toric affine surface. To describe its fan we need to introduce some notation. Let $N\simeq \mathbb{Z}^2$ be the lattice of 1-parameter subgroups of the torus $T:=\mathbb{C}^\ast\times \mathbb{C}^\ast$. Fix a $\mathbb{Z}$-basis $\{e_1,e_2\}$ of $N$ and define the vector $v_i:=ie_1+(1-i)e_2\in N$ for any integer $i\geq 0$. Then the fan of $\mathbb{C}^2/\mathbb{Z}_k$ consists of the two-dimensional cone $\sigma:=\operatorname{Cone}(v_0,v_k)\subset N_{\mathbb{Q}}$ and its subcones. The origin is the unique singular point of $\mathbb{C}^2/\mathbb{Z}_k$, and is a particular case of the so-called {\em rational double point} or {\em Du Val singularity} (see~\cite[Def.~10.4.10]{book:coxlittleschenck2011}).

By~\cite[Example~10.1.9 and Cor.~10.4.9]{book:coxlittleschenck2011}, the minimal resolution of singularities of $\mathbb{C}^2/\mathbb{Z}_k$ is the smooth toric surface $\varphi_k\colon X_k\to\mathbb{C}^2/\mathbb{Z}_k$ defined by the fan $\Sigma_k\subset N_{\mathbb{Q}}$, where
\begin{align}
\Sigma_k(0)&:=\big\{\{0\}\big\}\ ,\\[4pt]
\Sigma_k(1)&:=\big\{\rho_i:=\operatorname{Cone}(v_i)\,\big\vert\, i=0,1,2, \ldots,
k \big\}\ ,\\[4pt]
\Sigma_k(2)&:=\big\{\sigma_i:=\operatorname{Cone}(v_{i-1}, v_i)\,\big\vert\, i=1,2,
\ldots, k\big\} \ . 
\end{align}
Note that the vectors $v_i$ are the minimal generators of the rays $\rho_i$ for $i=0, 1, \ldots, k$. 

Let us consider the vector $b_\infty:=-v_0-v_k=-ke_1+(k-2)e_2$ in $N.$ Denote by $\rho_\infty$ the ray $\operatorname{Cone}(b_\infty)\subset N_\mathbb{R}$ and by $v_\infty$ its minimal generator. For even $k$, $v_\infty=\frac{1}{2}b_\infty$; for odd $k$, $v_\infty=b_\infty.$ Let $\sigma_{\infty,k}$ and $\sigma_{\infty,0}$ be the two-dimensional cones $\operatorname{Cone}(v_{k},v_\infty)\subset N_\mathbb{R}$ and $\operatorname{Cone}(v_{0}, v_\infty)\subset N_\mathbb{R}$ respectively. Let $\bar{X}_k$ be the normal projective toric surface defined by the fan $\bar{\Sigma}_k\subset N_\mathbb{R}$:
\begin{align}
\bar{\Sigma}_k(0)&:=\big\{\{0\} \big\} \ = \ \Sigma_k(0) \ ,\\[4pt]
\bar{\Sigma}_k(1)&:=\{\rho_i\,\vert\,i=0,1,2, \ldots,
k\}\cup\{\rho_\infty\} \ = \ \Sigma_k(1)\cup\{\rho_\infty\}\ ,\\[4pt]
\bar{\Sigma}_k(2)&:=\{\sigma_i\,\vert\, i=1,2, \ldots,
k\}\cup\{\sigma_{\infty,k},\sigma_{\infty,0}\} \ = \ \Sigma_k(2)\cup\{\sigma_{\infty,k}, \sigma_{\infty,0}\}\ .
\end{align}

First note that $i\colon X_k\hookrightarrow \bar{X}_k$ as an open dense subset. We denote by $D_\infty$ the $T$-invariant divisor associated to the ray $\rho_\infty$. 

\subsubsection*{Canonical orbifold}

Let $\tilde{k}$ be $k$ for odd $k$, otherwise $k/2$. Let $\pi^{\mathrm{can}}_k\colon \mathscr{X}_k^{\mathrm{can}}\to \bar{X}_k$ be the so-called {\em canonical toric orbifold} over $\bar{X}_k$ with torus $T$. It is the unique (up to isomorphism) smooth two-dimensional separated toric\footnote{Following \cite{art:fantechimannnironi2010}, a {\em Deligne-Mumford torus} $\mathscr{T}$ is a product $T\times \mathscr{B} G$ where $T$ is a ordinary torus and $G$ a finite abelian group. A {\em smooth toric Deligne-Mumford stack} is a smooth separated Deligne-Mumford algebraic stack $\mathscr{X}$ of finite type over $\mathbb{C}$, with as a coarse moduli space $\pi\colon \mathscr{X}\to X$ a scheme, together with an open immersion of a Deligne-Mumford torus $\imath\colon \mathscr{T}\hookrightarrow \mathscr{X}$ with dense image such that the action of $\mathscr{T}$ on itself extends to an action $a\colon \mathscr{T}\times \mathscr{X}\to \mathscr{X}$.} Deligne-Mumford algebraic stack of finite type over $\mathbb{C}$ such that the locus over which $\pi^{\mathrm{can}}_k$ is not an isomorphism has a nonpositive dimension.

As a global quotient stack, $\mathscr{X}_k^{\mathrm{can}}$ is isomorphic to $\left[Z_{\bar{\Sigma}_k}/(\mathbb{C}^\ast)^k\right]$, where $Z_{\bar{\Sigma}_k}$ is the union over all cones $\sigma\in\bar{\Sigma}_k$ of the open subsets
\begin{equation}\label{eq:Zsigma}
Z_\sigma:=\big\{x\in\mathbb{C}^{k+2}\, \big\vert \, x_i\neq0 \ \mbox{ if
} \ \rho_i\notin\sigma \big\} \ \subset \ \mathbb{C}^{k+2}\ ,
\end{equation}
and the action of $(\mathbb{C}^\ast)^k$ is given by
\begin{equation*}
(t_1,\dots,t_k)\cdot(z_1,\dots,z_{k+2})=\left\{
\begin{array}{ll}
\Big(\,\prod\limits_{i=1}^{k-1\, }t_i^i\, t_k^{2-k}\,
z_1\,,\,\prod\limits_{i=1}^{k-1}\, t_i^{-(i+1)}\, t_k^k\, z_2\,,\,t_1\,
z_3\,,\,\dots\,,\,t_k\, z_{k+2}\, \Big) \ , & \ k\mbox{ odd}\\[8pt]
\Big(\,\prod\limits_{i=1}^{k-1}\, t_i^i\, t_k^{1-\tilde{k}}\,
z_1\,,\,\prod\limits_{i=1}^{k-1}\, t_i^{-(i+1)}\, t_k^{\tilde{k}}\,
z_2\,,\,t_1\, z_3\,,\,\dots\,,\,t_k\, z_{k+2}\, \Big) \ , & \ k\mbox{ even}
\end{array}\right.
\end{equation*}
for  $(t_1,\dots,t_k)\in (\mathbb{C}^\ast)^k$ and $(z_1,\dots,z_{k+2})\in Z_{\bar{\Sigma}_k}$.

The effective Cartier divisor $\tilde{\mathscr{D}}_\infty:=(\pi_k^{\mathrm{can}})^{-1}(D_\infty)_{\mathrm{red}}$ is a toric orbifold with torus $\mathbb{C}^\ast$. By \cite[Prop.~3.10]{art:bruzzopedrinisalaszabo2013} $\tilde{\mathscr{D}}_\infty$ is isomorphic as a quotient stack to $[\mathbb{C}^2\setminus\{0\}/\mathbb{C}^\ast\times\mu_{\tilde{k}}]$, where the action is given by $(t,\omega)\cdot(z_1,z_2)=(t\omega\, z_1,t\, z_2)$ for $(t,\omega)\in\mathbb{C}^\ast\times\mu_{\tilde{k}}$ and $(z_1,z_2)\in\mathbb{C}^2\setminus\{0\}$. Moreover, $D_\infty\simeq \mathbb{P}^1$ is the coarse moduli space of $\tilde{\mathscr{D}}_\infty$ and the line bundle $\mathcal{O}_{\mathscr{X}_k^{\mathrm{can}}}(\tilde{\mathscr{D}}_\infty)$ is ample (cf.\ \cite[Rem.~3.18]{art:bruzzopedrinisalaszabo2013}).
\begin{remark}
As explained in \cite[Sect.~3.3]{art:bruzzopedrinisalaszabo2013}, $\tilde{\mathscr{D}}_\infty$ is an orbifold of type $(\tilde{k},\tilde{k})$ with two orbifold points $\tilde{p}_0$ and $\tilde{p}_\infty$, which are $(\pi_k^{\mathrm{can}})^{-1}(0)_{\mathrm{red}}$ and $(\pi_k^{\mathrm{can}})^{-1}(\infty)_{\mathrm{red}}$ respectively, where $0, \infty \in D_\infty$ are the two torus fixed points. By Proposition \ref{prop:fundamentalgroup-orbifoldcurves}, the fundamental group of $\tilde{\mathscr{D}}$ is $\mathbb{Z}_{\tilde{k}}$. Moreover, its universal cover is $\mathscr{P}(\frac{\tilde{k}}{\tilde{k}},\frac{\tilde{k}}{\tilde{k}})=\mathscr{P}(1,1)\simeq\mathbb{P}^1$.
\end{remark}

\subsubsection*{Root stack over the canonical orbifold}
 
Let $\phi_k\colon\mathscr{X}_k\to\mathscr{X}_k^{\mathrm{can}}$ be the stack obtained from $\mathscr{X}_k^{\mathrm{can}}$ by performing a $k$-\emph{root construction}\footnote{For the theory of root stacks we refer to \cite{art:cadman2007}.} along the divisor $\tilde{\mathscr{D}}_\infty$. As explained in \cite[Sect.~3.4]{art:bruzzopedrinisalaszabo2013}, $\mathscr{X}_k$ is a two-dimensional toric orbifold with torus $T$ and with coarse moduli space $\pi_k:=\pi_k^{\mathrm{can}}\circ\phi_k\colon \mathscr{X}_k\to\bar{X}_k$. As a global quotient stack, $\mathscr{X}_k$ is isomorphic to $[Z_{\bar{\Sigma}_k}/(\mathbb{C}^\ast)^k]$, where the action of $(\mathbb{C}^\ast)^k$ on $Z_{\bar{\Sigma}_k}$ is
\begin{equation*}
(t_1, \ldots, t_k)\cdot(z_1,\dots,z_{k+2})=\left\{\begin{array}{ll}
\Big(\, \prod\limits_{i=1}^{k-1}\, t_i^i\,
t_k^{2k-k^2}\,z_1,\,\prod\limits_{i=1}^{k-1}\, t_i^{-(i+1)}\,
t_k^{k^2}\,z_2,\,t_1\,z_3,\,\dots\,,\,t_k\,z_k\, \Big)&\mbox{ for odd } k\ ,\\[8pt]
\Big(\, \prod\limits_{i=1}^{k-1}\, t_i^i\, t_k^{k-k\,
  \tilde{k}}\,z_1,\,\prod\limits_{i=1}^{k-1}\, t_i^{-(i+1)}\, t_k^{k\,
  \tilde{k}}\,z_2,\,t_1\,z_3,\,\dots\,,\,t_k\,z_k\, \Big)&\mbox{ for even } k\ .
\end{array}\right.
\end{equation*}
for $(t_1,\dots,t_k)\in (\mathbb{C}^\ast)^k$ and $(z_1,\dots,z_{k+2})\in Z_{\bar{\Sigma}_k}$. 

Let $\mathscr{D}_\infty$ be the effective Cartier divisor $\pi_k^{-1}(D_\infty)_{\mathrm{red}}$; it is a smooth toric Deligne-Mumford stack with Deligne-Mumford torus $\mathscr{T}\simeq \mathbb{C}^\ast\times\mathscr{B}\mu_k$. By \cite[Prop.~3.30]{art:bruzzopedrinisalaszabo2013} $\mathscr{D}_\infty$ is isomorphic as a global quotient stack to
\begin{equation}\label{eq:quotientstack}
\left[\frac{\mathbb{C}^2\setminus\{0\}}{\mathbb{C}^\ast\times\mu_k}\right]\ ,
\end{equation}
where the action is given by
\begin{equation}\label{eq:action-infty}
(t,\omega)\cdot(z_1,z_2)=\left\{
\begin{array}{ll}
\big(t^{\tilde{k}}\, \omega \, z_1\,,\,t^{\tilde{k}}\, \omega^{-1}\,
z_2 \big) & \mbox{ for even $k$}\ ,\\
\big(t^k\, \omega^{\frac{k+1}{2}}\, z_1\,,\,t^k\,
\omega^{\frac{k-1}{2}}\, z_2\big) & \mbox{ for odd $k$}\ ,
\end{array}
\right.
\end{equation}
for $(t,\omega)\in\mathbb{C}^\ast\times\mu_k$ and $(z_1,z_2)\in\mathbb{C}^2\setminus\{0\}$. Moreover, $\mathscr{D}_\infty$ is a $\mu_k$-gerbe over $\tilde{\mathscr{D}}_\infty$.

Under this description of $\mathscr{D}_\infty$, the restriction of the line bundle $\mathcal{O}_{\mathscr{X}_k}(\mathscr{D}_\infty)$ on $\mathscr{D}_\infty$ is associated with the character of $\mathbb{C}^\ast\times\mu_k$ given by the projection to the first factor (cf.\ \cite[Lem.~3.35]{art:bruzzopedrinisalaszabo2013}). So it is easy to see that $\mathcal{O}_{\mathscr{X}_k}(\mathscr{D}_\infty)$ is ample on $\mathscr{D}_\infty$ (see also \cite[Sect.~6.1]{art:bruzzosala2013}). Moreover, (cf.\ \cite[Sect.~3.4]{art:bruzzopedrinisalaszabo2013})
\begin{equation}
\mathcal{O}_{\mathscr{X}_k}(\mathscr{D}_\infty)_{\vert \mathscr{D}_\infty}^{\otimes \tilde{k}k}\,\simeq \,\pi_k^\ast \mathcal{O}_{\bar{X}_k}(\tilde{k}D_\infty)_{\vert D_\infty}\,\simeq \, {\pi_k}_{\vert \mathscr{D}_\infty}^\ast \mathcal{O}_{D_\infty}(k/\tilde{k})\ .
\end{equation}
So by Proposition \ref{prop:positiveness} the line bundle $\mathcal{O}_{\mathscr{X}_k}(\mathscr{D}_\infty)$ is positive on $\mathscr{D}_\infty$.

\begin{remark} 
Since $\mathscr{D}_\infty^{\mathrm{top}}$ is a $\mu_k$-gerbe over $\tilde{\mathscr{D}}_\infty^{\mathrm{top}}$, 
its universal cover is $\mathscr{P}(\tilde{k},\tilde{k})$ by Proposition \ref{prop:fundamentalgroup-analyticcurves}. 
The group $\mathbb{C}^\ast\times\mu_k$ fits into a triangle as \eqref{eq:triangle} where the homomorphism 
$\rho\colon \mathbb{C}^\ast\times\mu_k\to \mathrm{GL}(2, \mathbb{C})$ is given by the action \eqref{eq:action-infty} and the homomorphism
$\mathbb{C}^\ast\to \mathbb{C}^\ast\times\mu_k$ is simply $t\mapsto (t, 1)$. Thus by the arguments in Section \ref{sec:fundamental-cover}, we obtain
that $\pi_1(\mathscr{D}_\infty^{\mathrm{top}})\simeq \mathbb{Z}_{k}$.
\end{remark}

\begin{remark}
 It is easily seen that $\mathscr{U}_{\mathscr{D}_{\infty}}$ and $\mathscr{U}_{\tilde{\mathscr{D}}_{\infty}}$ are actually diffeomorphic manifolds. Their homotopy exact sequences fit into the diagram:
 \begin{equation}
  \begin{tikzpicture}[xscale=2.2,yscale=-1.2]
    \node (A0_0) at (0.4, 0) {$1$};
    \node (A1_0) at (1, 0) {$C_{\mathscr{D}_{\infty}}$};
    \node (A2_0) at (2, 0) {$\pi_1(\mathscr{U}_{\mathscr{D}_{\infty}})$};
    \node (A3_0) at (3, 0) {$\pi_1(\mathscr{D}_{\infty}) $};
    \node (A4_0) at (3.7, 0) {$1$};

        \node (A0_1) at (0.4, 1) {$1$};
    \node (A1_1) at (1, 1) {$C_{\tilde{\mathscr{D}}_{\infty}}$};
    \node (A2_1) at (2, 1) {$\pi_1(\mathscr{U}_{\tilde{\mathscr{D}}_{\infty}})$};
    \node (A3_1) at (3, 1) {$\pi_1(\tilde{\mathscr{D}}_{\infty}) $};
    \node (A4_1) at (3.7, 1) {$1$};
    
    \path (A0_0) edge [->]node [auto] {$\scriptstyle{}$} (A1_0);
      \path (A1_0) edge [->]node [auto] {$\scriptstyle{}$} (A2_0);
           \path (A2_0) edge [->]node [auto] {$\scriptstyle{}$} (A3_0);
            \path (A3_0) edge [->]node [auto] {$\scriptstyle{}$} (A4_0);
                
  \path (A0_1) edge [->]node [auto] {$\scriptstyle{}$} (A1_1);
       \path (A1_1) edge [->]node [auto] {$\scriptstyle{}$} (A2_1);
           \path (A2_1) edge [->]node [auto] {$\scriptstyle{}$} (A3_1);
               \path (A3_1) edge [->]node [auto] {$\scriptstyle{}$} (A4_1);          
                
      \path (A1_0) edge [->]node [auto] {$\scriptstyle{}$} (A1_1);             
            \path (A2_0) edge [->]node [auto] {$\scriptstyle{\mathrm{id}}$} (A2_1);  
                  \path (A3_0) edge [->]node [auto] {$\scriptstyle{}$} (A3_1);  
  \end{tikzpicture} 
\end{equation} 
where the rightmost vertical arrow is a surjection and the leftmost is an injection. Since $\pi_1(\mathscr{D}_\infty^{\mathrm{top}})\simeq \mathbb{Z}_{k}$ and $\pi_1(\tilde{\mathscr{D}}_\infty^{\mathrm{top}})\simeq \mathbb{Z}_{\tilde k}$, then $C_{\tilde{\mathscr{D}}_{\infty}}\slash C_{\mathscr{D}_{\infty}}=\mathbb{Z}_{k/\tilde{k}}$. 
\end{remark}

By Theorem \ref{thm:bando} and Corollary \ref{cor:bando} we obtain the following.
\begin{theorem}\label{thm:ALE}
Given a holomorphic vector bundle $\mathcal E$ on $\mathscr{X}_k$, its restriction $\mathcal E_{\vert\mathscr{D}_\infty}$ is isomorphic to a fixed vector bundle $\mathcal F_\infty$ endowed with a flat unitary connection $\nabla$ iff there exists a Hermite-Einstein vector bundle $(\mathcal E', h')$ on ${ (X_k,\omega_k)}$ such that
$\displaystyle\int_{\mathscr{X}_k} |\Theta(\mathcal E',h')|^2 <+\infty$  and $(\mathcal{E}',h')$ has at infinity the holonomy 
given by the holonomy of $\nabla$, where $\omega_k$ is the ALE metric on $X_k$, and $\mathcal{E}$ is an extension of $\mathcal{E}'$.
\end{theorem}

\begin{rem}
The original result of Bando permits only to treat the case  where $k=2$ and the holonomy is trivial on 
$C_{\tilde{\mathscr{D}}_{\infty}}\slash C_{\mathscr{D}_{\infty}}=\mathbb{Z}_{2}$. \hfill $\triangle$
\end{rem}

\end{document}